\documentclass[a4paper,11pt]{article}

\usepackage{amsfonts}
\usepackage{amsmath}
\usepackage{amssymb}
\usepackage{graphicx}
\usepackage{amscd}
\usepackage{tikz}
\usepackage{epsfig}
\usepackage{psfrag}

\title{More on the properties of the generalized majorization}

\author{Marija Dodig\thanks{CEAFEL, Departamento de Mat\'ematica, Universidade de Lisboa, Edificio C6, Campo Grande, 
1749-016 Lisbon, Portugal, and
Mathematical Institute SANU, Knez Mihajlova 36, 11000 Belgrade, Serbia. ({\tt msdodig@fc.ul.pt}). Corresponding author.}
\and
Marko Sto\v si\'c\thanks{CAMGSD, Departamento de Matem\'atica, 
Instituto Superior T\'ecnico,
Av. Rovisco Pais 1, 1049-001 Lisbon, Portugal, and
Mathematical Institute SANU, Knez Mihajlova 36, 11000 Belgrade, Serbia. }
}

\date{}

\newtheorem{theorem}{Theorem}

\newtheorem{corollary}[theorem]{Corollary}

\newtheorem{definition}{Definition}

\newtheorem{lemma}{\indent Lemma}

\newtheorem{problem}{Problem}

\def\max{\mathop{\rm max}}

\def\kraj{\hfill\rule{6pt}{6pt}}

\begin{document}

\maketitle

\begin{abstract}
In this paper, we give corrected and  improved definitions of the sets $S$ and $\Delta$ compared to \cite{ela}. By using  these new definitions, we go throughout the proof of the main result in \cite{ela}, and we correct it. 
\end{abstract}

\section{Introduction}

\begin{definition}\label{majq}Let $d_1\ge \cdots\ge d_{m+k-s}$, $g_1\ge \cdots\ge g_{m+k}$, $a_1\ge \cdots\ge a_s$ be integers. Consider partitions $\mathbf{d}=(d_1,\ldots,d_{m+k-s}),$
 $\mathbf{g}=(g_1,\ldots,g_{m+k})$ and $\mathbf{a}=(a_1,\ldots,a_s).$ 
If \begin{eqnarray}
&d_i\ge g_{i+s},\quad\quad i=1,\ldots, m+k-s,\label{gm21}\\
&\sum_{i=1}^{h_j}{g_i}-\sum_{i=1}^{h_j-j}{d_i}\le \sum_{i=1}^j{a_i}, \quad\quad j=1,\ldots,s\label{gm11}\\
&\sum_{i=1}^{m+k}{g_i}=\sum_{i=1}^{m+k-s}{d_i}+\sum_{i=1}^s{a_i},\label{gm41}
\end{eqnarray}
where $$h_j:=\min\{i|d_{i-j+1}<g_i\},\quad j=1,\ldots,s,$$
 then we say that 
$\mathbf{g}$ 
is {\it majorized} by  $\mathbf{d}$ and $\mathbf{a}$. This type of majorization we  call {\it the generalized majorization}, and we write 
$$\mathbf{g}\prec'(\mathbf{d},\mathbf{a}).$$
\end{definition}

Notice that, if (\ref{gm41}) is satisfied, then (\ref{gm11}) is equivalent to the following:
\begin{equation}\sum_{i=h_j+1}^{m+k}{g_i}\ge \sum_{i=h_j-j+1}^{m+k-s}{d_i}+\sum_{i=j+1}^s{a_i},\quad j=1,\ldots,s.\label{gm31}
\end{equation}

\begin{definition}\label{weak}
If partitions $\mathbf{a}$, $\mathbf{d}$ and $\mathbf{g}$ in Definition \ref{majq} satisfy (\ref{gm21}), (\ref{gm31}) and $$\sum_{i=1}^{m+k}{g_i}\ge \sum_{i=1}^{m+k-s}{d_i}+\sum_{i=1}^s{a_i},$$ then we say that $\mathbf{g}$ 
is {\it weakly majorized} by  $\mathbf{d}$ and $\mathbf{a}$, and we write $$\mathbf{g}\prec''(\mathbf{d},\mathbf{a}).$$
\end{definition}
\begin{lemma}\label{zzz}(\cite[Lemma 2]{simax})
Suppose that  $\mathbf{d}=(d_1,\ldots, d_{m})$, $\mathbf{g}=(g_1,\ldots, g_{m+s})$ and $\mathbf{a}=(a_1,\ldots, a_s)$ satisfy $\mathbf{g}\prec''(\mathbf{d},\mathbf{a})$.  
Let $u$ be an integer such that $h_j< u  \le h_{j+1}$, for some $j\in\{0,\ldots,s\}$ ($h_0:=0$, $h_{s+1}:=m+s+1$). Then the following is also valid:
\begin{equation}\sum_{i=u}^{m+s}{g_i}\ge \sum_{i=u-j}^{m}{d_i}+\sum_{i=j+1}^s{a_i},\quad j=1,\ldots,s.\label{dokazi}\end{equation}
\end{lemma}

\begin{lemma}\label{2.4}\cite[Lemma 2.4]{ela}
 Let $\mathbf{a} = (a_1,\ldots,a_s)$, $\mathbf{d} = (d_1,\ldots,d_m)$, and $\mathbf{\bar{g}} = (\bar{g}_1,\ldots,\bar{g}_{m+s})$ be partitions such that 
 $$\mathbf{\bar{g}}\prec''(\mathbf{d},\mathbf{a}).$$
Let $f \in\{2, \ldots, m+s\}$, and let $\mathbf{g} = (g_1,
 \ldots, g_{m+s})$ be a partition such that
$$g_i=\bar{g}_i,\quad  i\ge f,$$
$$g_i\le \bar{g}_i,\quad  i<f,$$
$$\bar{g}_{f-1} \ge g_1 \ge g_{f- 1} \ge g_1 -1,$$
 $$\sum_{i=1}^{m+s}g_i \ge \sum_{i=1}^{m}d_i +
  \sum_{i=1}^{s}a_i.$$
Then $$\mathbf{{g}}\prec''(\mathbf{d},\mathbf{a}).$$

\end{lemma}

In \cite{ela} we have studied  the following problem:

\begin{problem}\label{p1} Let $m,n,s$ and $k$ be nonnegative integers such that $m+s=n+k$. Let $\mathbf{a}=(a_1,\ldots,a_s)$, $\mathbf{b}=(b_1,\ldots,b_k)$, $\mathbf{c}=(c_1,\ldots,c_n),$ and $\mathbf{d}=(d_1,\ldots,d_{m})$
be partitions.

Find necessary and sufficient conditions  for the existence of a partition $\mathbf{g}=(g_1,\ldots,g_{m+s})$, such that
$$\mathbf{g}\prec'(\mathbf{c},\mathbf{b})\quad\textrm{ and }\quad\mathbf{g}\prec'(\mathbf{d},\mathbf{a}).$$
\end{problem}

In fact, first we have resolved the following sub-problem: 

\begin{problem}\label{p12} Let $m,n,s$ and $k$ be nonnegative integers such that $m+s=n+k$. Let $\mathbf{a}=(a_1,\ldots,a_s)$, $\mathbf{b}=(b_1,\ldots,b_k)$, $\mathbf{c}=(c_1,\ldots,c_n),$ and $\mathbf{d}=(d_1,\ldots,d_{m})$
be partitions.

Find necessary and sufficient conditions  for the existence of a partition $\mathbf{g}=(g_1,\ldots,g_{m+s})$, such that
$$\mathbf{g}\prec''(\mathbf{c},\mathbf{b})\quad\textrm{ and }\quad\mathbf{g}\prec''(\mathbf{d},\mathbf{a}).$$
\end{problem} 


By Proposition 2.6 in \cite{ela} from now on we shall consider partitions $\mathbf{c}$ and $\mathbf{d}$ such that $c_i\ne d_j$ for all $i=1,\ldots,n$, and all $j=1,\ldots,m.$\\

Although we have solved Problem \ref{p1} in Theorem 5.1 from \cite{ela}, the solution strongly uses the definition of the sets $S$ and $\Delta$ from \cite{ela}, which  is not correct for all the values of $q_j$ and $q'_j$.  In this errata we are fixing all the problems in the definition of the sets  $S$ and $\Delta$ in \cite{ela}, and we give new, correct necessary and sufficient conditions for Problems \ref{p1} and \ref{p12}.

\section{Partitions and their properties}

Let $s,m,n$ and $k$ be positive integers such that 
$$m+s=n+k.$$
By a partition we assume a nonincreasing sequence of integers. In this paper we shall consider following partitions: 
\begin{eqnarray}
&\mathbf{a}=(a_1,\ldots,a_s)\label{a}\\
&\mathbf{d}=(d_1,\ldots,d_m)\label{d}\\
&\mathbf{b}=(b_1,\ldots,b_k)\label{b}\\
&\mathbf{c}=(c_1,\ldots,c_n),\label{c}
\end{eqnarray}
where $c_i\ne d_j$, for all $i=1,\ldots,n$
and $j=1,\ldots,m$. 



Denote by $\mathbf{u}$  the union of partitions $\mathbf{c}$ and $\mathbf{d}$, by $\mathbf{e}$  the union of partitions $\mathbf{d}$ and $\mathbf{a}$, and by  $\mathbf{e'}$  the union of partitions $\mathbf{c}$ and $\mathbf{b}$. Thus, we have

 $$\mathbf{u}=(u_1,\ldots,u_{n+m}):=(d_1,\ldots, d_m)\cup (c_1,\ldots,c_n),$$
$$\mathbf{e}=(e_1,\ldots,e_{m+s}):=(d_1,\ldots, d_m)\cup (a_1,\ldots,a_s),$$
and 
$$\mathbf{e'}=(e'_1,\ldots,e'_{m+s}):=(c_1,\ldots, c_n)\cup (b_1,\ldots,b_k).$$

In the definition of $e_i$'s, if $d_i=a_j$, then let $i_j=\min\{i|d_i=a_j\}$, and let $u=\min\{i|a_i=a_j\}, $ and $v=\max\{i|a_i=a_j\}$. Then we put $e_{i_j+u-1}=a_u$, $e_{i_j+u}=a_{u+1}$, \ldots, $e_{i_j+v}=a_v$, $e_{i_j+v+1}=d_{i_j}$ {{(i.e. $\mathbf{e}: \cdots a_u\ge\cdots\ge a_v\ge d_{i_j}\ge\cdots$)}}. Analogously, if $c_i=b_j$, then let $i_j=\min\{i|c_i=b_j\}$, and let $u=\min\{i|b_i=b_j\},$ and $v=\max\{i|b_i=b_j\}$. Then we put $e'_{i_j+u-1}=b_u$, $e'_{i_j+u}=b_{u+1}$, \ldots, $e'_{i_j+v}=b_v$, $e'_{i_j+v+1}=c_{i_j}$.

{{For any sequence of integers $y_1,\ldots,y_w$ we put $\sum_{i=r}^s y_i=0$ if $r>s$.}}
Moreover for any such sequence, we assume $y_i=+\infty$, for $i\le 0$, and $y_i=-\infty$, for $i>w$.
\subsection{New, improved definition of the sets $\mathbf{S}$ and $\mathbf{\Delta}$} 

In this section we improve the definition of the sets $S$ and $\Delta$ given in \cite{ela}. This is the main feature of this errata. After introducing these new and improved definitions, we are left with adjusting the main result in \cite{ela}, which will be done in the sequel sections.

\begin{definition}\label{sd}
Definition of the sets $S$ and $\Delta$ is given inductively. We start by putting $S$ and $\Delta$ to be empty sets, and then we fill them in the following way, step by step: \\

\noindent We start by choosing the smallest element in $\mathbf{u}$. If there are equals among $c_i$'s or $d_i$'s, we always first choose the  element with the largest index (note that we are assuming $c_i\ne d_j$ for all $i=1,\ldots,n$ and $j=1,\ldots,m$).\\

\paragraph{--} If the chosen element belongs to $\mathbf{d}$, say $d_j$, then we calculate
\begin{equation}q_j:= s-\sharp\{i\in S|c_i<d_j\}+\sharp\{i>j|i\notin \Delta\}+1.\label{qj}\end{equation}

\noindent Next we check the following:
\begin{eqnarray}
\bullet&& \textrm{If  }\quad q_j>s \quad\Rightarrow\quad\textrm{ then } j\in \Delta\nonumber\\
\bullet&&\textrm{If }\quad q_j\le s \quad\Rightarrow\quad\textrm{then let $l\in S$ be the minimal index such that $d_j>c_l$}\nonumber\\
&& (a)\quad \textrm{ Now, if  }\nonumber\\
&&\quad\quad\sharp\{i| a_i>c_l\}\ge s-\sharp\{i\in S|i>l\}+\sharp\{ i\notin \Delta|d_i<c_l\},\label{onoprvo}\\
&&\textrm{ and if  $d_j$ belongs to the smallest }\nonumber\\
&&\quad\quad \sharp\{i| a_i>c_l\}- s+\sharp\{i\in S|i>l\}-\sharp\{ i\notin \Delta|d_i<c_l\}+1\label{onoprvo1}\\
&&\textrm{  $e_i$'s bigger than $c_l$, then we put  $j\notin \Delta,$}\nonumber\\
&& (b) \quad \textrm {otherwise we check the inequality }\nonumber
\end{eqnarray}
\begin{equation}\label{q}\sum_{c_i<d_j, i\in S}c_i\ge \sum_{i\notin \Delta, i>j}d_i+d_j+\sum_{i=q_j+1}^sa_i.\end{equation}

If the equation (\ref{q}) is satisfied, then we put $j\notin \Delta$, and if the equation (\ref{q}) is not satisfied then we put $j\in \Delta$.\\

\paragraph{--} \quad If the chosen  element belongs to $\mathbf{c}$,
say $c_j$, then we have the dual definition, i.e.  we consider
\begin{equation}q'_j:= k-\sharp\{i\in \Delta|d_i<c_j\}+\sharp\{i>j|i\notin S\}+1.\label{qj'}\end{equation}

Then we check the following:
\begin{eqnarray}
\bullet&& \textrm{If  }\quad q'_j>k \quad\Rightarrow\quad\textrm{ then } j\in S\nonumber\\
\bullet&&\textrm{If }\quad q'_j\le k \quad\Rightarrow\quad\textrm{then let $l\in \Delta$ be the minimal index such that $c_i>d_l$}\nonumber\nonumber\\
&&\quad\quad  (a)\quad \textrm{ Now, if  }\nonumber\\
&&\quad\quad\sharp\{i| b_i>d_l\}\ge k-\sharp\{i\in \Delta|i>l\}+\sharp\{ i\notin S|c_i<d_l\},\label{onodrugo}\\
&&\quad\quad\textrm{ and if  $c_j$ belongs to the smallest }\nonumber\\
&&\quad\quad \sharp\{i| b_i>d_l\}- k+\sharp\{i\in \Delta|i>l\}-\sharp\{ i\notin S|c_i<d_l\}+1\label{onodrugo1}\\
&&\quad\quad   \textrm{   $e'_i$'s bigger than } d_l, \textrm{  then we put  } j\notin S\nonumber\\
&& \quad\quad (b) \quad \textrm {otherwise we check the inequality }\nonumber
\end{eqnarray}

\begin{equation}\label{q'}\sum_{d_i<c_j, i\in \Delta}d_i\ge \sum_{i\notin S, i>j}c_i+c_j+\sum_{i=q'_j+1}^kb_i.\end{equation}

If the equation (\ref{q'}) is satisfied, then we put $j\notin S$, and if the equation (\ref{q'}) is not satisfied then we put $j\in S$.\\


Now choose the next smallest element in $\mathbf{u}$, and proceed until  all the elements in $\mathbf{u}$ are checked.  This ends our definition of the sets $S$ and $\Delta$.\end{definition}\kraj


We note here, that the difference between Definition \ref{sd} and  the definition of the sets $S$ and $\Delta$ from \cite{ela}, is in indices $i$ and $j$ for which $q_i>s$ and $q'_j>k$. Also, there is improvement in the definition for the indices for which $q_i\le s$ and $q'_j\le k$ if (\ref{onoprvo}) and (\ref{onoprvo1}), and respectively, (\ref{onodrugo}) and (\ref{onodrugo1}) are valid.\\

Now, as in \cite{ela}, we re-name all  $d_i$'s with $i\in \Delta$, and call them $d^1\ge\cdots\ge d^h$, where $h=\sharp \Delta$. Analogously,  re-name all  $c_i$'s with $i\in S$, and call them $c^1\ge\cdots\ge c^{h'}$, where $h'=\sharp S$.

Analogously as in \cite{ela}, in order to simplify the notation, we define the following integers related to the sets $S$ and $\Delta$:

\begin{definition}\label{four}
For every $d^j$, $j=1,\ldots,h$, we define
$$m'_j:=\sharp\{i| b_i>d^j\} $$
$$t'_j:=k-(h-j)+\sharp\{i\notin S|c_i<d^j\}$$
$$z'_j:=\sharp\{i|c_i>d^j\},$$
and for every 
$c^j$, $j=1,\ldots,h'$, we define
$$m_j:=\sharp\{i| a_i>c^j\}$$
$$t_j:=s-(h'-j)+\sharp\{i\notin \Delta|d_i<c^j\}$$
$$z_j:=\sharp\{i|d_i>c^j\}.$$
\end{definition}
In addition, we also formally define $d^0:=d_0=+\infty,$ $d^{h+1}:=-\infty,$ $t'_{h+1}=k+1,$ $z'_{h+1}=n,$
 and we extend definitions of $m'_j$, $t'_j$ and $z'_j$ to the case $j=0$:
$m'_0:=\sharp\{i| b_i>d^0\}=0,$
$t'_0:=k-h+\sharp\{i\notin S|c_i<d^0\}=k-h+\sharp\{i\notin S\}=n+k-h-h',$ and 
$z'_0:=\sharp\{i|c_i>d^0\}=0.$

Analogously, we also formally define 
$c^0:=c_0=+\infty,$
$c^{h'+1}:=-\infty,$ $t_{h'+1}=s+1,$ $z_{h'+1}=m,$
and we extend definitions of $m_j$, $t_j$ and $z_j$ to the case $j=0$:
$m_0:=\sharp\{i| a_i>c^0\}=0,$
$t_0:=s-h'+\sharp\{i\notin \Delta|d_i<c^0\}=s-h'+\sharp\{i\notin \Delta\}=m+s-h-h',$
$z_0:=\sharp\{i|d_i>c^0\}=0.$

Note that since $m+s=n+k$, we have $t_0=t'_0$. Also, by Definition \ref{sd} we have $t'_{h}=k$ and $t_{h'}=s$.

\begin{definition}\label{w}
For $y\in\{0,\ldots,h'\}$ we define:
$$w_y:=\sharp\{i\notin \Delta| c^y>d_i>c^{y+1}\}.$$
For $x\in\{0,\ldots,h\}$ we define:
$$w'_x:=\sharp\{j\notin S| d^x>c_j>d^{x+1}\}.$$
\end{definition}
From Definitions \ref{four} and \ref{w} we directly obtain:
\begin{lemma}\label{dva}\begin{eqnarray}t_{x+1}&=&t_x+1-w_x,\quad  x=0,\ldots,h',\\ 
t'_{y+1}&=&t'_y+1-w'_y,\quad  y=0,\ldots,h,\\ 
z_x+t_x&<&z_{x+1}+t_{x+1},\quad  x=0,\ldots,h',\\
z'_y+t'_y&<&z'_{y+1}+t'_{y+1},\quad  y=0,\ldots,h.\end{eqnarray}
\end{lemma}

Now we can re-write the conditions (\ref{onoprvo}), (\ref{q}), (\ref{onodrugo}) and (\ref{q'}) in Definition \ref{sd} in the following way:\\

For $d_j$, $j\in\{1,\ldots,m\}$, let $l\in\{0,\ldots,h'\}$ be such that $c^{l}>d_j>c^{l+1}$. Then
$$q_j= s-(h'-l)+\sharp\{i>j|i\notin \Delta\}+1,$$
and condition (\ref{onoprvo}) becomes
$$m_{l+1}\ge t_{l+1},$$
and (\ref{q}) is equal to






\begin{equation}\label{quj}\sum_{i=l+1}^{h'}c^i\ge \sum_{i\notin \Delta, i>j}d_i+d_j+\sum_{i=q_j+1}^sa_i.\end{equation}



Analogously, for $c_j$, $j\in\{1,\ldots,n\}$, let 
 ${{l'}}\in\{0,\ldots,h\}$ be  such that $d^{{l'}}>c_j>d^{{l'+1}}$. Then
$$q'_j= k-(h-{{l'}})+\sharp\{i>j|i\notin S\}+1.$$
Also, (\ref{onodrugo}) becomes
$$m'_{{l'+1}}\ge t'_{{l'+1}},$$
and (\ref{q'}) is equal to


\begin{equation}\label{q'uj}\sum_{i={{l'+1}}}^{h}d^i\ge \sum_{i\notin S, i>j}c_i+c_j+\sum_{i=q'_j+1}^kb_i.\end{equation}





\section{Auxiliary lemmas}

In the following section we give auxiliary lemmas which are used  in the proof of the main result. In fact, some of these lemmas coincide with  lemmas from \cite{ela}. However, since we have changed  definition of the sets $S$ and $\Delta$, we have to prove them again. This is done for Lemmas 4.1, 4.3, {{4.5 and 4.6.}} Also, Lemmas 4.7 and 4.8 in \cite{ela} are now included in the definition of the sets $S$ and $\Delta$, while Lemmas 4.9 and 4.10 are included in Lemma \ref{dva}. The rest of the lemmas in \cite{ela} are not correct or necessary anymore.\\

In the rest of the paper we shall use the notation from Problem \ref{p1} and from Definitions \ref{sd}, \ref{four} and \ref{w}.
 
\begin{lemma}\label{lema1}\cite[Lemma 4.1]{ela}{ {Let $y\in\{0,\ldots,h'\}$ and let $j\in\{1,\ldots,m-1\}$ be such that }}  $c^y>d_j\ge d_{j+1}>c^{y+1}$.  Then, if $j+1\in \Delta$ we have that $j\in \Delta$.
\end{lemma}
\textbf{Proof:}
Since $j+1\in \Delta$, we have $q_j=q_{j+1}$. From the definition of $\Delta$, there are two possibilities: either $q_{j+1}>s$, and  then $q_j>s$, i.e.  $j\in \Delta$, as wanted; either (\ref{q}) is not valid for $d_{j+1}$, in which case we trivially obtain that it is not valid for $d_j$ as well. Hence
 $j\in \Delta$, as wanted.
\kraj

Completely analogously we have the dual result:
\begin{lemma}\label{lema2}\cite[Lemma 4.3]{ela}{ {Let $x\in\{0,\ldots,h\}$ and let $j\in\{1,\ldots,n-1\}$   be such that }}
 $d^x>c_j\ge c_{j+1}>d^{x+1}$.  Then if $j+1\in S$ we have that $j\in S$.
\end{lemma}
\begin{lemma}\cite[Lemma 4.6]{ela}\label{ttt} Let $j \in \Delta$. Let $i \in \{1,\ldots,h\}$ be such that $d_j = d^i$ and let $x \in \{0,\ldots,h'\}$ be such that $c^x > d_j > c^{x+1}$. Then
$$z'_i + t'_i = j + t_x.$$\end{lemma}
\textbf{Proof:} By Definition \ref{four}, together with Lemmas \ref{lema1} and \ref{lema2}, we obtain 
$$z'_i+t'_i=\sharp\{l|c_l>d^i\}+k-(h-i)+\sharp\{l\notin S|c_l<d^i\}=$$
$$=
k-(h'-i)+(n-\sharp\{l\in S|c_l<d^i\})=
k-(h-i)+n-(h'-x)=$$
$$=m+s-(h-i)-(h'-x)=s-(h'-x)+(m-\sharp\{l\in \Delta|l>j\})=$$
$$=
s-(h'-x)+j+\sharp\{l\notin \Delta|l>j\}=
j+s-(h'-x)+\sharp\{l\notin \Delta|c^x>d_l\}=j+t_x.$$
\kraj\\
Dually, we have
\begin{lemma}\cite[Lemma 4.5]{ela} \label{zaz2}Let $j\in S$. Let $i\in\{1,\ldots,h'\}$ be such that $c_j =c^i$ and let $x\in 
\{0,\ldots,h\}$ be such that $d^x > c_j > d^{x+1}$. Then
$$z_i + t_i = j + t'_x.$$
\end{lemma}

To proceed we also need the following lemma from \cite{simax}:
\begin{lemma}\label{najnajnov}[Lemma 4.9 \cite{simax}]
Let $u_1\ge\cdots\ge u_k$ and $v_1\ge\cdots\ge v_k$ be integers. If 
$$\sharp \{i|u_i> v_j\}\ge j,\quad \textrm{ for all } \quad j=1,\ldots,k,$$
then 
$$\sum_{i=1}^k{u_i}\ge \sum_{i=1}^k{v_i}+k.$$ 
\end{lemma}

{{\begin{lemma}\label{tnula}Let $j\in\{1,\ldots,m\}$ be such that $j\in \Delta$. Let  $y\in\{0,\ldots,h'\}$ be such that $c^y>d_j>c^{y+1}.$ Then $t_y\ge 0$.
\end{lemma}
\textbf{Proof:}
 Indeed, if $t_y<0$ then :
$$m_{y+1}-t_{y+1}+1=m_{y+1}-t_y-1+w_y+1>m_{y+1}+w_y.$$
The last means that $d_{z_{y+1}-w_y}$ is among the smallest $m_{y+1}-t_{y+1}+1$ $\quad e_i$'s larger than $c^{y+1}$. Since, by Lemma \ref{u1}, we have that  $q_{z_{y+1}-w_y}\le s$, by the part $(a)$ of the definition of the set $\Delta$, we conclude $z_{y+1}-w_y\notin \Delta$, which is a contradiction by the definition of $w_y$. Hence $t_y\ge 0$, as wanted. \kraj\\

Dually, we have 
\begin{lemma}\label{t'nula}Let $j\in\{1,\ldots,n\}$ be such that $j\in S$. Let  $x\in\{0,\ldots,h\}$ be such that $d^x>c_j>d^{x+1}.$ Then $t'_x\ge 0$.
\end{lemma}
}}

\begin{lemma}\label{t0}
$t_0=t'_0\ge 0$.
\end{lemma}
\textbf{Proof:}
If any of the sets $S$ or $\Delta$ is empty, we directly get that $t_0\ge0$. 
If none of the sets $S$ and $\Delta$ is empty,  we have that if $d^1> c^1$ by Lemma \ref{tnula} $t_0\ge 0$, and if $c^1> d^1$ by Lemma \ref{t'nula} $t'_0\ge 0$, as wanted.
\kraj\\

Lemmas \ref{ttt}, \ref{zaz2} and \ref{t0} together give:
\begin{lemma}\label{lemnov}
The numbers $z_{i}+t_i$ for $i=1,\ldots,h'$, and $z'_{i}+t'_i$ for $i=1,\ldots,h$, are all distinct. In addition, $$\{z_i+t_i| i=1,\ldots,h'\}\cup\{z'_i+t'_i| i=1,\ldots,h\}=\{t_0+1,t_0+2,\ldots,m+s\}.$$\end{lemma}

\subsection{Novel lemmas}

Next, we give two new lemmas comparing to \cite{ela}. They will play important role in the main result: 

\begin{lemma}\label{u1}
Suppose that $c^{h'}\ge a_s$, and let {{$j\in\{1,\ldots,m\}$}} be such that $d_j>c^{h'}$. Then $q_j\le s$. {{In addition if $j\notin \Delta$ then $q_j<s$.}}

\end{lemma}

\textbf{Proof: }   Before proceeding note that by the definition of $q_l$ all $d_l<c^{h'}$ satisfy $l\in \Delta$.

Since $d_j>c^{h'}$, we have that $1\le j\le z_{h'}$. Let $p\in\{0,\ldots,h'-1\}$ be  such that $c^{p}>d_{j}>c^{p+1}$.  The rest of the proof goes by the induction on $j$.\\

Let  $j=z_{h'}$. By  definition  (\ref{qj}), we have $q_{z_{h'}}=s-(h'-p)+1\le s$, as wanted.\\

Now let $1\le  j <  z_{h'}$ and suppose that  $q_i\le s$, for all $i=j+1,\ldots,z_{h'}$. We shall prove that then $q_{j}\le s$.

{{
By definition (\ref{qj}), we have that if $q_{j+1}<s$, then $q_{j}\le s$. So the only case we are left to consider is  when $q_{j+1}=s$.

 Let $y\in\{0,\ldots,h'-1\}$ be  such that $c^{y}>d_{j+1}>c^{y+1}$, and  let  
$$\gamma=\sharp\{i\notin \Delta| i=j+2,\ldots,z_{y+1}\}.$$

 We shall prove that $j+1\in \Delta$, and then by definition (\ref{qj}) will follow 
$$q_{j}\le q_{j+1}= s, \textrm{ as wanted}.$$

Since $t_{y+1}=q_{j+1}-\gamma=s-\gamma$, and $m_{y+1}\le s-1$ (since $c^{h'}\ge a_s$), we have  $m_{y+1}-t_{y+1}+1\le \gamma$, so by the definition of $\gamma$ we have that $d_{j+1}$ doesn't satisfy part $(a)$ of the definition of the set $\Delta$. So we are left with checking the condition $(b)$ of the definition of the set $\Delta$, i.e. we are left with checking 
\begin{equation}\sum_{i=y+1}^{h'} c^i<\sum_{i=j+2, i\notin \Delta}^m d_i+d_{j+1}.\label{uf}\end{equation}



Let  $h'-y=1+\sharp\{i\notin \Delta | j+2\le i\le m\}$ (since $q_{j+1}=s$). Let $u_1\ge\cdots\ge u_{h'-y}$ be the non increasing ordering of $d_{j+1}$ and $d_i$ with $j+1\le i\le m$, $i\notin \Delta$, and let $v_1\ge\cdots\ge v_{h'-y}$ be defined as $v_i:=c^{y+i}$, $i=1,\ldots,h'-y$. We claim that  $u_i> v_i$, $i=1,\ldots,h'-y$.\\

Since $d_{j+1}<c^{y+1}$ we have $u_1>v_1$. Now let us fix $i_0\in\{2,\ldots,h'-y\}$. Then $u_{i_0}=d_l$ for some $l\notin \Delta$ with $j+2\le l\le m, $ i.e. $i_0=1+\sharp\{i\notin\Delta| j+2\le i\le l\}$. Let $r\in\{0,\ldots,h'-1\}$ be such that $c^r>d_l>c^{r+1}$. Note that $l\le z_{h'}$ since for all $i>z_{h'}$ we have $i\in\Delta$.

From $q_l\le s$ we get
\begin{equation}\label{zv1}
\sharp\{i\ge l | i\notin\Delta\}\le h'-r.\end{equation}
On the other hand, $q_{j+1}=s$ gives
\begin{equation}\label{zv2}
1+\sharp\{i\notin \Delta | i\ge j+2\}= h'-y.\end{equation}
Then (\ref{zv1}) and (\ref{zv2}) together give
$$
1+\sharp\{i\notin \Delta | j+2<i\le l\}\ge r+1-y,$$
i.e. $$i_0\ge r+1-y.$$
Therefore
$$u_{i_0}=d_l>c^{r+1}=c^{y+(r+1-y)}\ge c^{y+i_0}=v_{i_0},$$ as wanted. 
Then by Lemma \ref{najnajnov} we get (\ref{uf}). Thus, we have proved that $j+1\in \Delta$, and so $q_{j}\le q_{j+1}= s$, as wanted.}}

\kraj
 
Dually, we get :
\begin{lemma}\label{u2}
Suppose that $d^{h}\ge b_k$, and let {{$j\in\{1,\ldots,n\}$}} be such that $c_j>d^{h}$. Then $q'_j\le k$. {{In addition if $j\notin S$ then $q'_j<k$.}}

\end{lemma}

\vskip 0.3cm


As direct corollaries of Lemmas \ref{u1} and \ref{u2}, we have 
\begin{corollary}\label{cor1}\begin{equation}\label{txs}c^{h'}\ge a_s \quad \Longrightarrow \quad t_y<s,\quad \textrm{ for all }\quad y=0,\ldots,h'-1,\end{equation}
\begin{equation}\label{tpxs}d^{h}\ge b_k \quad \Longrightarrow \quad  t'_x<k, \quad \textrm{ for all }\quad  x=0,\ldots,h-1.\end{equation}
\end{corollary}
\textbf{Proof:}
We shall prove (\ref{txs}), and (\ref{tpxs}) follows dually. 

First note that there are no $i\notin\Delta$ such that $c^{h'-1}>d_i>c^{h'}$. Indeed, suppose on the contrary that  $j\in\{1,\ldots,m\}$ is the largest such index. Since $m_{h'}\le s-1$ and $t_{h'}=s$, $j\notin\Delta$ implies that (\ref{q}) is satisfied, i.e.
$c^{h'}\ge d_j$ which is a contradiction. Therefore $t_{h'-1}=s-1$. 

Now fix $y\in\{0,\ldots,h'-2\}$. If there are no $i\notin\Delta$ such that $c^y>d_i>c^{h'-1}$ then $t_y=t_{h'-1}-(h'-1-y)=s-1-(h'-1-y)<s$. If there exists $i\notin \Delta$ with $c^y>d_i>c^{h'-1}$, then let $j$ be the smallest such index and let $p\in\{y,\ldots,h'-2\}$ be such that $c^p>d_j>c^{p+1}$. Then $t_p=q_j$, and so by Lemma \ref{u1} $t_y=t_p-(p-y)=q_j-(p-y)<s-(p-y)\le s$, as wanted.

\kraj

\subsection{A partition mutually generally majorized by two pairs of partitions}

Consider the partitions $\mathbf{a,d,b}$ and $\mathbf{c}$ as in (\ref{a})--(\ref{c}). In this subsection we shall assume that there exists a  partition $\mathbf{g}=(g_1,\ldots,g_{m+s})$, such that 
\begin{equation}\mathbf{g}\prec''(\mathbf{d},\mathbf{a})\quad\textrm{ and }\quad 
\mathbf{g}\prec''(\mathbf{c},\mathbf{b}).\label{est11u}\end{equation}
Under this assumption, we prove the following four lemmas (all together they correct and prove analogous results to Lemmas 5.2, 5.3, 5.4 and 5.5. from \cite{ela}):
\begin{lemma}\label{L12}Let $\mathbf{a,d,b,c}$ and $\mathbf{g}$  be partitions which satisfy
 (\ref{est11u}). Then
\begin{equation} c^{h'}\ge g_{z_{h'}+s} \quad \textrm{ and }  \quad d^{h}\ge g_{z'_{h}+k},  \label{dokazati}\end{equation}\label{dokazatiL}
as well as \begin{equation} c^{h'}\ge a_s \quad \textrm{ and }  \quad d^{h}\ge b_k.  \label{dokazatiu}\end{equation}
\end{lemma}
\textbf{Proof:}  We shall prove that $ c^{h'}\ge g_{z_{h'}+s}$ and $c^{h'}\ge a_s$, and the proof of $ d^{h}\ge g_{z'_{h}+k}$ and $d^{h}\ge b_k$ goes completely dually, by changing the roles of the partitions $\mathbf{c}$ and  $\mathbf{d}$, as well as  $\mathbf{a}$ and  $\mathbf{b}$, respectively.\\

If suppose that  $d_m> c^{h'}$, i.e. if $z_{h'}=m$,  then $c^{h'}=c_n$ and since 
$\mathbf{g}\prec''(\mathbf{c},\mathbf{b})$
we have 
$$c^{h'}=c_n\ge g_{n+k}=g_{m+s}=g_{z_{h'}+s}, \textrm{ as wanted.}$$

If $z_{h'}<m$, then $c^{h'}=c_{n-\alpha+1}$ for some $1\le \alpha\le n$, and $z_{h'}=m-\beta$, for some $1\le\beta\le m.$ Then we have that $i\notin S$ for $n-\alpha+1<i\le n$, and $j\in \Delta$ for $m-\beta<j\le m$.\\

If  $\beta<  \alpha$, we have $c^{h'}=c_{n-\alpha+1}\ge g_{n-\alpha+1+k}=g_{m-\alpha+1+s}\ge g_{m-\beta+s}=g_{z_{h'}+s} $, as wanted.\\

If $\beta\ge \alpha$, then 
from the definition of $q'_i$ we have $$q'_{n-\alpha+1}=k-\beta+\alpha\le k.$$
 Since $n-\alpha+1\in S$,  from the definition of the set $S$ (part $(a)$) we have that the index $n-\alpha+1$ does not belong to the $m'_{h-\beta+1}-t'_{h-\beta+1}+1$ smallest $e'_i$'s bigger than $d_{m-\beta+1} (=d_{z_{h'}+1})$. Let 
 $$\bar{u}=\sharp \{i\in\{1,\ldots,k\}| b_i> c_{n-\alpha+1}\}, $$
$$\bar{v}=\sharp \{i\in\{1,\ldots,k\}| c_{n-\alpha+1}\ge b_i> d_{m-\beta+1}\}, $$
$$\bar{w} =\sharp \{n-\alpha+1<i\le n| c_i> d_{m-\beta+1}\} $$
and $$\bar{z}=\sharp \{n-\alpha+1<i\le n| c_i <d_{m-\beta+1}\} .$$

Then $\bar{z}+\bar{w}=\alpha-1$, $t'_{h-\beta+1}=k-(\beta-1)+\bar{z}$ and $m'_{h-\beta+1}=\bar{u}+\bar{v}$. Since $n-\alpha+1\in S$ we have $\bar{v}+\bar{w}\ge m'_{h-\beta+1}-t'_{h-\beta+1}+1=\bar{u}+\bar{v}-k+\beta-\bar{z}$, i.e. $\bar{u}\le \bar{w}+\bar{z}+k-\beta=\alpha-1+k-\beta$. Thus, 
$$\alpha+k>\beta,$$
and 
\begin{equation}c^{h'}=c_{n-\alpha+1}\ge b_{\alpha+k-\beta}\label{2}\end{equation}
Also, since $n-\alpha+1\in S$ by the part $(b)$ of the definition of the set $S $ (since $q'_{n-\alpha+1}\le k$), we have 
\begin{equation}\label{pp1}\sum_{i=m-\beta+1}^m d_i< c_{n-\alpha+1}+\sum_{i=n-\alpha+2}^n c_i+\sum_{i=k+\alpha-\beta+1}^k b_i.\end{equation}



Now, let us  suppose the opposite from what we need to prove, i.e. that $c^{h'}< g_{z_{h'}+s}.$
Last is equivalent to $c_{n-\alpha+1}< g_{m-\beta+s}$. Thus, by definition of $h'_j=\min\{i|c_{i-j+1}<g_i\}$, we have  $h'_{m-\beta+s-n+\alpha}\le m+s-\beta$, i.e. $h'_{k+\alpha-\beta}\le m+s-\beta$. Let $u\in\{0,\ldots,k\}$ be such that $h'_u\le m+s-\beta < h'_{u+1}$. Then $u\ge k+\alpha-\beta$.\\

Since $\mathbf{g}\prec''(\mathbf{c},\mathbf{b})$,  by the definition of the weak generalized majorization, and by Lemma \ref{zzz}, we have 
\begin{equation}\label{avg}\sum_{i=m+s-\beta+1}^{m+s}g_i\ge \sum_{i=m+s-\beta+1-u}^n c_i +\sum_{i=u+1}^k b_i.\end{equation} 

Since  $\mathbf{g}\prec''(\mathbf{d},\mathbf{a})$ implies $d_i\ge g_{i+s}$, $i=1,\ldots,m$, by (\ref{avg}) we have
\begin{equation}\sum_{i=m-\beta+1}^{m}d_i\ge \sum_{i=m+s-\beta+1-u}^n c_i +\sum_{i=u+1}^k b_i.\label{3}\end{equation}




Since  $u\ge k+\alpha-\beta$, from (\ref{2}) we have that 
$$\sum_{i=m+s-\beta+1-u}^n c_i+\sum_{i=u+1}^k b_i=\sum_{i=n-\alpha+1}^n c_i+\sum_{i=k-\alpha+\beta+1}^k b_i$$
$$ +(\sum_{i=m+s-\beta+1-u}^{n-\alpha} c_i-\sum_{i=k+\alpha-\beta+1}^u b_i)\ge$$
$$\ge \sum_{i=n-\alpha+1}^n c_i+\sum_{i=k-\alpha+\beta+1}^k b_i,$$ which together with (\ref{3}) gives
\begin{equation}\label{pp2}\sum_{i=m-\beta+1}^m d_i\ge \sum_{i=n-\alpha+1}^n c_i+\sum_{i=k+\alpha-\beta+1}^k b_i, 
\end{equation}
which contradicts (\ref{pp1}). Thus,  $c^{h'}\ge g_{z_{h'}+s}$.\\

Now, let us prove that $c^{h'}\ge a_s$. Let $j\in\{0,\ldots,s\}$, be  such that $h_j<z_{h'}+s\le h_{j+1}$ ($h_0=0$, $h_{s+1}=m+s+1$). 
Then $\mathbf{g}\prec''(\mathbf{d},\mathbf{a})$  (by Lemma \ref{zzz} and the definition of the weak generalized majorization) gives
\begin{equation}\sum_{i=z_{h'}+s}^{m+s}g_i\ge \sum_{i=z_{h'}+s-j}^m d_i+\sum_{i=j+1}^s a_i.\label{ast}\end{equation}

Equations (\ref{gm21}) and (\ref{ast}) together with $c^{h'}\ge g_{z_{h'}+s}$ give

\begin{equation}\label{nov}c^{h'}+\sum_{i=z_{h'}+1}^m d_i\ge \sum_{i=z_{h'}+s-j}^m d_i+\sum_{i=j+1}^s a_i.\end{equation}
If $j=s$, (\ref{nov}) becomes $c^{h'}\ge d_{z_{h'}}$ which is a contradiction by the definition of $z_{h'}$. On the other hand if $j<s$, then  (\ref{nov}) gives

$$(s-j)c^{h'}\ge c^{h'}+\sum_{i=z_{h'}+1}^{z_{h'}+s-j-1} d_i\ge \sum_{i=j+1}^s a_i\ge (s-j) a_s ,$$
i.e. $c^{h'}\ge a_s$, as wanted.

\kraj\\

\begin{lemma}\label{vaz1} Let $\mathbf{a,d,b,c}$ and $\mathbf{g}$  be partitions which satisfy $\mathbf{g}\prec''(\mathbf{d},\mathbf{a})$.
Suppose that  $c^{h'}\ge a_s$. Let $j\in\{1,\ldots,m\}$ be such that $j\in \Delta$. Let  $y\in\{0,\ldots,h'\}$ be such that $c^y>d_j>c^{y+1}.$\\

If $$c^l\ge g_{z_l+t_l},\quad  \textrm{ for all }\quad l\ge y+1,$$
 and 
 $$d_{\alpha}\ge g_{\alpha+t_{\beta}}, \quad \textrm{ for all }\quad \alpha\in \Delta,\quad  \alpha>z_{y+1},\textrm{ 
 and }
c^{\beta}>d_{\alpha}>c^{\beta+1},$$
 then 
\begin{equation}\label{222}
d_j\ge g_{j+t_y}.
\end{equation}

\end{lemma}

\textbf{Proof:}
If $y=h'$, we have that $t_y=t_{h'}=s$, and so (\ref{222}) becomes $d_j\ge g_{j+s}$, which follows from $\mathbf{g}\prec''(\mathbf{d},\mathbf{a})$.\\

So, from now on, we assume $0\le y \le h'-1$. Since $c^{h'}\ge a_s$, by (\ref{txs}) we have $t_y<s$.  Also, by Lemma \ref{tnula} we  have that $t_y\ge 0$. 
Therefore, we have $0\le t_y <s$.
We shall prove that 
\begin{equation}h_{t_y+1}\ge z_{y+1}+t_{y+1},\label{4}\end{equation} 
where $h_{t_y+1}=\min\{u|d_{u-t_y}<g_u\}$. 
If (\ref{4}) is valid then $d_u\ge g_{u+t_y}$, for $u+t_y< z_{y+1}+t_{y+1}$, i.e. $u\le z_{y+1}+t_{y+1}-t_y-1=z_{y+1}-w_y$, thus proving  the lemma.\\

Let suppose the opposite to (\ref{4}), i.e. let $h_{t_y+1}\le z_{y+1}+t_{y+1}-1$. Let $u\in\{1,\ldots,s\}$ be such that $h_u<z_{y+1}+t_{y+1}\le h_{u+1}$. Then $u\ge t_y+1$ and since $\mathbf{g}\prec''(\mathbf{d},\mathbf{a})$, by the definition of the weak generalized majorization, and by Lemma \ref{zzz}, we have:
 \begin{equation}\sum_{i=z_{y+1}+t_{y+1}}^{m+s}g_i\ge\sum_{i=z_{y+1}+t_{y+1}-u}^md_i+\sum_{i=u+1}^s a_i.\label{5uj}\end{equation}

By the assumptions of the lemma, we have
\begin{equation}\sum_{i=y+1}^{h'}c^i+\sum_{j\in \Delta, j>z_{y+1}}d_j\ge  \sum_{i=z_{y+1}+t_{y+1}}^{m+s}g_i.\label{5uj1}\end{equation}

Inequalities (\ref{5uj}) and (\ref{5uj1}), together give 
\begin{equation}\sum_{i=y+1}^{h'}c^i+\sum_{j\in \Delta, j>z_{y+1}}d_j\ge \sum_{i=z_{y+1}+t_{y+1}-u}^md_i+\sum_{i=u+1}^s a_i.\label{5}\end{equation}


Since $z_{y+1}-w_y\in \Delta$, and since $q_{z_{y+1}-w_y}=t_y+1\le s$ , 
we have that  $d_{z_{y+1}-w_y}$  does not satisfy the condition  from the part $(b)$ of the definition of the set $\Delta$:
$$\sum_{i=y+1}^{h'} c^i<d_{z_{y+1}-w_y}+\sum_{i>z_{y+1}-w_y, i\notin \Delta}d_i+\sum_{i=t_y+2}^s a_i$$
 which further gives
$$\sum_{i=y+1}^{h'} c^i+\sum_{i>z_{y+1}, i\in \Delta}d_i <\sum_{i=z_{y+1}-w_y}^m d_i+\sum_{i=t_y+2}^s a_i$$
Last equation together with (\ref{5}) give
$$\sum_{i=z_{y+1}+t_{y+1}-u}^m d_i+\sum_{i=u+1}^s a_i<\sum_{i=z_{y+1}-w_y}^m d_i+\sum_{i=t_y+2}^s a_i.$$
Since $u\ge t_y+1$ and $t_y=t_{y+1}-1+w_y$, we have
\begin{equation}\label{contr1}
\sum_{i=z_{y+1}+t_{y+1}-u}^{z_{y+1}-w_y-1}d_i<\sum_{i=t_y+2}^{u} a_i.
\end{equation}
Note that there is the same number of summands on the left and the right hand side in (\ref{contr1}). 
Since $z_{y+1}-w_y\in \Delta$, we know that $d_{z_{y+1}-w_y}$ does not belong to the smallest $m_{y+1}-t_{y+1}+1$ $e_i$'s larger than $c^{y+1}$. Therefore  
$m_{y+1}-t_{y+1}+1\le w_y+\sharp\{i| d_{z_{y+1}-w_y}>a_i> c^{y+1}\},$  i.e. $\sharp\{i|a_i\ge d_{z_{y+1}-w_y}\}\le t_y$. This is equivalent to $d_{z_{y+1}-w_y}> a_{t_y+1}$, and so the smallest summand on the LHS of (\ref{contr1}) is larger then the largest summand on the RHS,  which gives a contradiction. Thus (\ref{4}) is valid, and so we have proved our lemma.

\kraj\\

Dually, we have:

\begin{lemma}\label{vaz2} Consider partitions $\mathbf{a,b,g,d},$ and $\mathbf{c}$. 
Let $\mathbf{g}\prec''(\mathbf{c},\mathbf{b})$. Suppose that $d^{h}\ge b_k$.
 Let $j\in\{1,\ldots,n\}$ be such that $j\in S$. Let  $x\in\{0,\ldots,h\}$ be such that $d^x>c_j>d^{x+1}.$

If $$d^l\ge c_{z'_l+t'_l}, \quad \textrm{ for all }\quad l\ge x+1,$$ and 
$$c_{\alpha}\ge g_{\alpha+t'_{\beta}}, \textrm{ for all } \alpha\in S, \quad \alpha>z'_{x+1},\textrm{ and } d^{\beta}>c_{\alpha}>d^{\beta+1},$$ 
then
\begin{equation}\label{222n}
c_j\ge g_{j+t'_x}.
\end{equation}
\end{lemma}

Next, we shall unify results from Lemmas \ref{L12} -- \ref{vaz2} and proving that if there exists a partition $\mathbf{g}$ satisfying $\mathbf{g}\prec''(\mathbf{d},\mathbf{a})$ and $\mathbf{g}\prec''(\mathbf{c},\mathbf{b}),$ that then $g_i$'s are bounded above by $c_i$'s with $i\in S$ and $d_j$'s with $j\in \Delta$. More precisely, we have:
\begin{lemma}\label{cetvrta} Let $\mathbf{a,d,b,c}$ and $\mathbf{g}$  be partitions which satisfy
 (\ref{est11u}). Then
\begin{eqnarray}
\label{nej1}c^i\ge g_{z_i+t_i},\quad i=1,\ldots,h',\\
\label{nej2}d^i\ge g_{z'_i+t'_i},\quad i=1,\ldots,h.
\end{eqnarray} 

\end{lemma}

\textbf{Proof:} 
Before proceeding, by (\ref{est11u}) and by Lemma \ref{L12} we have that $c^{h'}\ge a_s$ and $d^{h}\ge b_k.$ Thus, we can apply Lemmas \ref{vaz1} and \ref{vaz2}.

Next,  we note that  (\ref{nej2}) can be written in the following (equivalent) way: 

Since $d^i$ corresponds to  $d_{j}$  for some $j\in\{1,\ldots,m\}$, (i.e. $d^i=d_j$), let $y\in\{0,\ldots,h'\}$ be such that $c^y>d_j>c^{y+1}$. Then by Lemma \ref{ttt} (\ref{nej2}) can be equivalently written as 
\begin{equation}\label{nej3}
d_j\ge g_{j+t_y}.
\end{equation}

We can rewrite (\ref{nej1}) analogously:  if $c^i$ corresponds to  $c_{j}$ (i.e. $c^i=c_j$), for some $j\in\{1,\ldots,n\}$,  let $x\in\{0,\ldots,h\}$ be such that $d^x>c_j>d^{x+1}$. Then (\ref{nej1}) can be equivalently written as 
\begin{equation}\label{nej4}
c_j\ge g_{j+t'_x}.
\end{equation}

We shall prove  inequalities (\ref{nej1}) and (\ref{nej2}) together and by induction. More precisely, let $A$ be the union of $\{c^i|i=1,\ldots,h'\}$ and $\{d^i|i=1,\ldots,h\}$. Then the goal is to prove that each element of $A$ is larger or equal than certain $g_l$, for appropriate index $l$ in accordance with (\ref{nej1}) and (\ref{nej2}). We shall prove these inequalities by induction on the elements of $A$ by starting from the smallest element of $A$. In the process we observe the equal elements of $A$ in the order determined by the indices of $c^i$ and $d^i$, i.e. if for some $i$ we have $c^i=c^{i+1}$ we shall first prove it for $c^{i+1}$ and then for $c^i$ (recall that we are assuming that there no $i$ and $j$ with $c_i=d_j$).\\

Now, the base of induction is to prove the inequalities (\ref{nej1}) and (\ref{nej2}) for the smallest element of $A$, i.e. (\ref{nej1}) for $c^{h'}$, in the case $c^{h'}<d^{h}$, and (\ref{nej2}) for $d^{h}$, in the case $c^{h'}>d^{h}$.

If  $c^{h'}<d^{h}$, we have that  $c^{h'}=c_n$, $z_{h'}=m$ and $t_{h'}=s$, 
and (\ref{nej1}) becomes 
$$c_n\ge g_{n+k},$$
which follows by
$\mathbf{g}\prec''(\mathbf{c},\mathbf{b})$.

Analogously, if $c^{h'}>d^{h}$,  we have that in fact $d^{h}=d_m$, $z'_{h}=n$ and $t'_{h}=k$, 
and  (\ref{nej2}) becomes 
$$d_m\ge g_{m+s},$$
which follows by
$\mathbf{g}\prec''(\mathbf{d},\mathbf{a})$.\\



The induction step is proved in Lemmas \ref{vaz1} and \ref{vaz2}. Lemma \ref{vaz1} solves the case when the element from $A$ is $d^i$ for some $i\in\{1,\ldots,h\}$, and it proves that (\ref{nej2}) is valid for that $d^i$, if the inequalities (\ref{nej1}) and (\ref{nej2}) hold for all elements of $A$ smaller than $d^i$.
 
Lemma \ref{vaz2} solves the case when the element from $A$ is $c^i$ for some $i\in\{1,\ldots,h'\}$, and it proves that (\ref{nej1}) is valid for that $c^i$, if the inequalities (\ref{nej1}) and (\ref{nej2}) hold for all elements of $A$ smaller than $c^i$. 

Therefore, together with the above base of induction, Lemmas \ref{vaz1} and \ref{vaz2}, prove the inequalities (\ref{nej1}) and (\ref{nej2}).

\kraj

\section{Main result}

Now we can give our main result. It is very similar to the result in \cite{ela}, but here we cover all the possible cases, some of which were missing in \cite{ela}.\\

The following theorem resolves Problem \ref{p12}:
\begin{theorem}\label{glavna}
Let $\mathbf{a}$, $\mathbf{d}$, $\mathbf{b}$ and $\mathbf{c}$ be partitions as in (\ref{a})--(\ref{c}).  There exists a partition $\mathbf{g}=(g_1,\ldots,g_{m+s})$, such that 
\begin{equation}\label{prvi}\mathbf{g}\prec''(\mathbf{d},\mathbf{a})\quad\textrm{ and }\quad\mathbf{g}\prec''(\mathbf{c},\mathbf{b})\end{equation}
if and only if the following conditions are valid
\begin{eqnarray*}
(i)&&\textrm{if } \quad y\in\{1,\ldots, h'\}\quad  \textrm { is such that }\quad  t_y\le m_y \quad \textrm{ then}\\
&&\sum_{i=z_y+t_y}^{z_y+m_y}e_i\le \sum_{i=y}^{h'} c^i-\sum_{i\ge z_y+1, i\notin \Delta }d_i-\sum_{i=m_y+1}^s a_i,\\
(ii)&&\textrm{if } \quad x\in\{1,\ldots, h\}\quad  \textrm { is such that }\quad  t'_x\le m'_x \quad \textrm{ then}\\
&&\sum_{i=z'_x+t'_x}^{z'_x+m'_x}e'_i\le \sum_{i=x}^{h} d^i-\sum_{i\ge z'_x+1,i\notin S }c_i-\sum_{i=m'_x+1}^k b_i.
\end{eqnarray*}
\end{theorem}

\vskip 0.5cm

The following theorem resolves Problem \ref{p1}:

\begin{theorem}\label{glavna1}
Let $\mathbf{a}$, $\mathbf{d}$, $\mathbf{b}$ and $\mathbf{c}$ be partitions as in (\ref{a})--(\ref{c}).
  There exists a partition $\mathbf{g}=(g_1,\ldots,g_{m+s})$, such that 
\begin{equation}\mathbf{g}\prec'(\mathbf{d},\mathbf{a})\quad\textrm{ and }\quad\mathbf{g}\prec'(\mathbf{c},\mathbf{b})\label{drugi}\end{equation}
if and only if  \begin{equation}\sum_{i=1}^n{c_i}+\sum_{i=1}^k{b_i}=\sum_{i=1}^m{d_i}+\sum_{i=1}^s{a_i},\label{sumasuma}\end{equation} and the  conditions $(i)$ and $(ii)$ are valid.



\end{theorem}

Proofs of Theorems \ref{glavna} and \ref{glavna1} are given in the sequel sections. 


\section{Proof of Theorem \ref{glavna}}

\subsection{Necessity of conditions $(i)$ and $(ii)$}\label{nes}
Let  us  assume that there exists a partition $\mathbf{g}$ such that 
\begin{equation}\mathbf{g}\prec''(\mathbf{d},\mathbf{a})\label{est}\end{equation}
\begin{equation}\mathbf{g}\prec''(\mathbf{c},\mathbf{b}).\label{est11}\end{equation}
Then we shall prove  that  conditions  $(i)$ and $(ii)$ hold.\\

Before proceeding, we note that for all $j$ such that $c^{h'}>d_j$, we have $q_j>s$ and thus $j\in \Delta$. So we have 
$$c^{h'}>d_{z_{h'}+1}\ge\cdots\ge d_m  \Rightarrow z_{h'}+1,\ldots,m\in \Delta.$$

Also, for all $j$ such that $d^{h}>c_j$, we have $q'_j>k$ and thus $j\in S$. So we have 
$$d^{h}>c_{z'_{h}+1}\ge\cdots\ge c_n  \Rightarrow z'_{h}+1,\ldots,n\in S.$$

Let $y\in\{1,\ldots,h'\}$ be such that $t_y\le m_y.$ Let $u\in\{0,\ldots,s\}$ be such that $h_u<z_y+t_y\le h_{u+1}$ ($h_0=0$, $h_{s+1}=m+s+1$). From $\mathbf{g}\prec''(\mathbf{d},\mathbf{a})$, by the definition of the weak generalized majorization, and by Lemma \ref{zzz}, we have 
$$\sum_{i=z_y+t_y}^{m+s}g_i\ge \sum_{i=z_y+t_y-u}^m d_i+\sum_{i=u+1}^s a_i$$
Together with Lemma \ref{cetvrta} this gives
\begin{equation}\sum_{i=y}^{h'} c^i+\sum_{i>z_y, i\in \Delta}d_i\ge \sum_{i=z_y+t_y-u}^m d_i+\sum_{i=u+1}^sa_i.\label{9}\end{equation}
We need to consider three cases:
\begin{eqnarray}
&u<t_y\le m_y\label{6}\\
&t_y\le u< m_y\label{7}\\
&t_y\le m_y\le u\label{8}
\end{eqnarray}
 
For each of the cases we can write (\ref{9}) in the following form (for all details see the proof of formula (5.26) in \cite{ela}):

$$\sum_{i=y}^{h'} c^i-\sum_{i>z_y, i\notin \Delta}d_i -\sum_{i=m_y+1}^s a_i\ge \sum_{i=z_y+t_y}^{z_y+m_y}e_i,$$

which is exactly the condition $(i)$.\\



Completely analogously, by changing roles of $\mathbf{c}$ and $\mathbf{b}$ with $\mathbf{d}$ and $\mathbf{a}$, respectively, we obtain the dual result, i.e. we prove condition $(ii)$. This finishes the proof of the necessity of conditions.

\subsection{Sufficiency of conditions $(i)$ and $(ii)$}\label{suf}

Suppose now that  conditions $(i)$ and $(ii)$ are valid. 
In this section we shall define a partition $\mathbf{{g}}$ which satisfies
\begin{eqnarray}
\mathbf{{g}}\prec''(\mathbf{d},\mathbf{a})\quad\textrm{ and }\quad
\mathbf{{g}}\prec''(\mathbf{c},\mathbf{b}).\label{16q}
\end{eqnarray}








Before proceeding, we shall prove that conditions $(i)$ and $(ii)$ imply
\begin{equation}c^{h'}\ge a_s, \quad d^{h}\ge b_k,\label{novnov}\end{equation}

\begin{equation} \sum_{i=1}^{h'} c^i\ge \sum_{i\notin \Delta}d_i+\sum_{i=t_0+1}^s a_i.\label{312} \end{equation}

and

\begin{equation} \sum_{i=1}^{h} d^i\ge \sum_{i\notin S}c_i+\sum_{i=t'_0+1}^k b_i.\label{313} \end{equation}

First note that inequality $c^{h'}\ge a_s$ is equivalent  to $m_{h'}<t_{h'}=s,$ and inequality $d^{h}\ge b_k$ is equivalent to $m'_h<t'_h=k$. 

Suppose on the contrary that $s\le m_{h'}$, i.e. $m_{h'}=s$. Then by condition $(i)$ for $y=h'$ we would have
$$c^{h'}\ge e_{z_{h'}+s}=e_{z_{h'}+m_{h'}},$$ which is a contradiction. Analogously if $m'_h=k$
by condition $(ii)$ for $x=h$ we would have
$$d^{h}\ge e'_{z'_{h}+k},$$ which is a contradiction. Therefore 
$c^{h'}\ge a_s$ and $d^{h}\ge b_k$, as wanted.\\

Next, we shall prove (\ref{312}) -- the inequality (\ref{313}) is obtained completely dually.\\

Let $(i)$ and $(ii)$ be valid. 

First we suppose that there are no $i\in \{1,\ldots,m\}$ such that  $i\notin \Delta$. Then by the definition we have $t_0=s-h'$ and $t_i=t_{i-1}+1=t_0+i$, $i=1,\ldots,h'$. 

{{
If $m_i<t_i$ for all $i\in\{1,\ldots,h'\}$, then by the definition of $m_i$ we have $c^i\ge a_{t_i}=a_{t_0+i}$, and thus 
$$\sum_{i=1}^{h'} c^i\ge \sum_{i=t_0+1}^s a_i,$$
which is precisely (\ref{312}) in this case. }}


If there is $i\in\{1,\ldots,h'\}$ for which $m_i\ge t_i$, then let  $y\in\{1,\ldots,h'\}$ be the minimal such index.
Then condition $(i)$ for $c^y$ gives
\begin{equation}\label{mart}\sum_{i=z_y+t_y}^{z_y+m_y}e_i\le \sum_{i=y}^{h'} c^i-\sum_{i=m_y+1}^s a_i. \end{equation}
Among $e_i$'s on the LHS there can be no $d_i$, since by the part $(a)$ of the definition of the set $\Delta$, we would have that those $i\notin \Delta$, contradicting the assumption that there are no such indices. Therefore those $e_i$'s are precisely $a_{t_y},\ldots,a_{m_y}$  (note that $t_y=t_0+y\ge y>0$, by condition $(i)$), and so (\ref{mart}) is equal to
\begin{equation}\label{mart1}\sum_{i=y}^{h'} c^i\ge \sum_{i=t_y}^s a_i=\sum_{i=t_0+y}^s a_i.\end{equation}

Since for all $i=1,\ldots,y-1$ we have $m_i+1\le t_i=t_0+i$, from the definition of $m_i$, we have
$c^i\ge a_{t_0+i}$, for $i=1,\ldots,y-1$. This together with  (\ref{mart1}) prove (\ref{312}) in this case.\\

Now suppose that there exists  $i\in\{1,\ldots,m\}$ such that $i\notin \Delta$. Let $j$ be the minimal such index. {{
By the definition of the set $\Delta$, we have that $q_j\le s$, and thus, by the definition of $q_j$, we conclude that $S$ is nonempty.

Since all $d_i<c^{h'}$ satisfy $i\in \Delta$, there exists  $y\in\{1,\ldots,h'\}$  such that 
$$c^{y-1}>d_j>c^y.$$
Then by the definition of $j$ we have $j=z_y-w_{y-1}+1$.
Also,  we have that $t_i=t_0+i$, for $i=1,\ldots,y-1$.

If there exists $i\in\{1,\ldots,y-1\}$, such that $m_i\ge t_i$, then denote by $x$ the minimal such index. Then in exactly the same way as in the first case (since there are no $i\notin \Delta$ with $d_i>c^{y-1}$), we obtain that  condition $(i)$ for $c^x$ implies
$$\sum_{i=x}^{h'} c^i\ge \sum_{i\notin \Delta}d_i+\sum_{i=t_x}^s a_i=\sum_{i\notin \Delta}d_i+\sum_{i=t_0+x}^s a_i.$$



\noindent Together with $c^i\ge a_{t_0+i}$, for $i=1,\ldots,x-1$, this proves (\ref{312}).

Thus, suppose that  $m_i<t_i$, for all $i=1,\ldots,y-1$, and therefore 
\begin{equation}\label{pom99}
c^i\ge a_{t_0+i},\quad\quad i=1,\ldots,y-1.
\end{equation}
Now, since $j\notin \Delta$, we have two possibilities from the definition of $\Delta$. If the part $(a)$ of the definition is satisfied, $d_j$ is among the smallest $m_y-t_y+1$ $\,\,e_i$'s larger than $c^y$. Thus, $j,j+1,\ldots,{z_y}\notin \Delta$, as well as  $t_y\le m_y$.

Then  condition $(i)$ for $c^y$ gives:
\begin{equation}\label{pom88}
\sum_{i=z_y+t_y}^{z_y+m_y}e_i\le \sum_{i=y}^{h'} c^i-\sum_{i>z_y,\, i\notin \Delta}d_i-\sum_{i=m_y+1}^s a_i. 
\end{equation}
By the above assumptions
 $(e_{z_y+t_y},\ldots,e_{z_y+m_y})$ consists of $w_{y-1}$ $\, d_i$'s, while the remaining $m_y-t_y+1-w_{y-1}=m_y-t_{y-1}$ are $a_i$'s, i.e. they are precisely $a_{t_{y-1}+1},\ldots,a_{m_y}$ (they are all larger than $c^y$). So, (\ref{pom88}) becomes:
\begin{equation}\label{pom77} 
\sum_{i=y}^h c^i \ge \sum_{i\notin \Delta}d_i+\sum_{i=t_{y-1}+1}^s a_i=\sum_{i\notin \Delta}d_i+\sum_{i=t_{0}+y}^s a_i. 
\end{equation}

On the other hand, if $j\notin \Delta$ because of the part $(b)$ of the definition of $\Delta$, then  \begin{equation}\label{pom66} 
\sum_{i=y}^{h'} c^i \ge \sum_{i\notin \Delta}d_i+\sum_{i=q_{j}+1}^s a_i. 
\end{equation}
Since from the definition of $q_i$'s and $t_i$'s we have that $q_j=t_{y-1}$, the last inequality becomes precisely (\ref{pom77}). }}

Therefore, we have obtained that (\ref{pom77}) holds, and together with (\ref{pom99}) finally gives the wanted condition (\ref{312}).     \\

 Completely analogously by changing the roles of partitions $\mathbf{c}$ and $\mathbf{b}$ with $\mathbf{d}$ and $\mathbf{a}$, respectively,  we obtain (\ref{313}).

\vskip 0.4cm

\subsubsection{Definition of $\mathbf{{g}}$}\label{barg}

Let $M=\max(a_1,b_1,c_1,d_1)+1$. By Lemma \ref{t0}, we have  $t_0=
m+s-(h+h')
\ge 0$.  Let $\mathbf{{g}}=({g}_1,\ldots,{g}_{m+s})$  be a partition defined as the following union 
$$\{c_i|i\in S\}\cup \{d_i|i\in \Delta\}\cup {\{M,\ldots,M\}}_{t_0}.$$


In other words we have
\begin{eqnarray}
\label{defbg1}{g}_1=\cdots={g}_{t_0}&=&\max(a_1,b_1,c_1,d_1)+1\\
\label{defbg3}{g}_{j}&=&d_{j-t_x},  \textrm{for } z_x+t_x<j<z_{x+1}+t_{x+1},  x=0,\ldots,h', \\
\label{defbg2}{g}_{z_x+t_x}&=&c^x,\quad x=1,\ldots,h'.
\end{eqnarray}
Equivalently we can write this also as
\begin{eqnarray}
{g}_1=\cdots={g}_{t'_0}&=&\max(a_1,b_1,c_1,d_1)+1\\
{g}_{j}&=&c_{j-t'_x}, \textrm{ for } z'_x+t'_x<j<z'_{x+1}+t'_{x+1}, x=0,\ldots,h,\\
{g}_{z'_x+t'_x}&=&d^x,\quad x=1,\ldots,h.
\end{eqnarray}

We shall prove that  $\mathbf{{g}}$ satisfies 
\begin{eqnarray}
\mathbf{{g}}&\prec''&(\mathbf{d},\mathbf{a})\label{15}\\
\mathbf{{g}}&\prec''&(\mathbf{c},\mathbf{b}).\label{16}
\end{eqnarray}

We start with proving (\ref{15}). By  Definition \ref{weak} of the weak majorization we need to prove the following:
\begin{eqnarray}
&d_i\ge {g}_{i+s}, \quad i=1,\ldots,m,\label{maj11}\\
&\sum_{i={h}_j+1}^{m+s} {g}_i \ge \sum_{i={h}_j-j+1}^m d_i + \sum_{i=j+1}^s a_i,\quad j=1,\ldots,s,\label{maj22}\\
&\sum_{i=1}^{m+s} {g}_i \ge \sum_{i=1}^m d_i +\sum_{i=1}^s a_i,\label{maj33}
\end{eqnarray}
where ${h}_j:=\min \{i|d_{i-j+1}<{g}_i\}$, for $j=1,\ldots,s$.\\

Regarding (\ref{maj11}),  since (\ref{novnov}) and (\ref{txs}) give  $t_0{{\le}} s$, we have that  ${g}_i$'s appearing in (\ref{maj11}) are the ones
defined by (\ref{defbg2}) and (\ref{defbg3}).

Now, if $i\in \Delta$, from (\ref{defbg3}) we have that $d_i={g}_{i+t_x}$, for some $x\in\{0,\ldots,h'\}$, and since {{$t_x\le s$
for any such $x$ we obtain $d_i\ge {g}_{i+s}$, as wanted. 

If on the other hand $i\notin \Delta$, then let $y\in\{0,\ldots,h'-1\}$ be such that $c^y>d_i>c^{y+1}$. Then we have that $i\in\{z_{y+1}-w_y+1,\ldots,z_y\}$, and by (\ref{defbg2}) we have:
$$d_i>c^{y+1}={g}_{z_{y+1}+t_{y+1}}={g}_{z_{y+1}-w_y+1+t_y}\ge {g}_{i+s},$$
since $z_{y+1}-w_y+1\le i$ and $t_y{{\le}} s$. This proves (\ref{maj11}).\\}}

Now, we pass to (\ref{maj22}).
First we note that from the definition of ${g}_i$, (\ref{defbg1})--(\ref{defbg3}), we can compute the values of ${h}_j$, for $j=1,\ldots,s$. We have that:
\begin{eqnarray}
&{h}_j=j,\quad j=1,\ldots,t_0,\label{barh1}\\
&{h}_j=z_x+t_x,\textrm{ where } x=\min\{i\in\{1,\ldots,h'\}|t_i=j\}, j=t_0+1,\ldots,s. \label{barh2}
\end{eqnarray}
Indeed, from (\ref{defbg1}) we have ${g}_{t_0}\ge d_1$, which gives (\ref{barh1}).\\

As for (\ref{barh2}) first note that $x$ is well-defined, i.e. the set $\{i\in\{1,\ldots,h'\}|t_i=j\}$ is non-empty, for $j=t_0+1,\ldots,s$. Indeed, from the definition of $t_x$, we have that $t_{x+1}=t_x+1-w_x$, and so $t_{x+1}\le t_x +1$, for $x=0,\ldots,h'-1$. Since $t_{h'}=s$, and $t_0\le s$ we have that the set $\{t_i|i=1,\ldots,h'\}$ contains all integers from the set $\{t_0+1,\ldots,s\}$.



Now, we show that for every $j
\in\{t_0+1,\ldots,s\}$, there exists  $i\in\{1,\ldots,h'\}$, such that ${{h}_j}=z_i+t_i$.


Indeed, if, on the contrary, there exists $j\in\{t_0+1,\ldots,s\}$, for which there are no   $i\in\{1,\ldots,h'\}$, such that ${{h}_j}=z_i+t_i$, then  let $u\in\{0,\ldots,h'\}$ be 
such that $z_u+t_u<{h}_j<z_{u+1}+t_{u+1}$. Then by (\ref{defbg3}) we have ${g}_{{h}_j}=d_{{h}_j-t_u}$, and from 
the definition of ${h}_j$, we have $d_{{h}_j-j+1}<{g}_{{h}_j}=d_{{h}_j-t_u}$, which implies $j\le t_u$,
and so $u\ge 1$. But then, from (\ref{defbg2}), ${g}_{z_u+t_u}=c^u>d_{z_u+1}\ge d_{z_u+t_u-j+1}$, and so 
${h}_j \le z_u+t_u$, which is a contradiction.

Hence we have that there exists $i\in\{1,\ldots,h'\}$ such that
${{h}_j}=z_i+t_i$. Then from the definition of ${h}_j$ we have $d_{z_i}>c^i={g}_{z_i+t_i}={g}_{{h}_j}>d_{{h}_j-j+1}=d_{z_i+t_i-j+1}$, and so $t_i\ge j$. Now, if $t_i>j$, since $t_{x+1}\le t_x +1$, for $x=0,\ldots,h'-1$, we have that there exists $u\in\{1,\ldots,i-1\}$ such that $t_u=j$. Then ${g}_{z_u+t_u}=c^u>d_{z_u+1}=d_{z_u+t_u-j+1}$, which together with $z_u+t_u<z_i+t_i$ (since $u<i$) contradicts the definition of  ${h}_j$. Therefore $t_i=j$ which finally proves (\ref{barh2}).\\

Now we shall prove (\ref{maj22}).

 Let  $j=1,\ldots,t_0$. By (\ref{barh1}), condition (\ref{maj22})  becomes  
 \begin{equation}\label{sept1}
\sum_{i=j+1}^{m+s}{g}_i\ge \sum_{i=1}^md_i+\sum_{i=j+1}^s a_i, \quad j=1,\ldots,t_0.
\end{equation}
By (\ref{defbg1}), it is enough to prove (\ref{sept1}) for $j=t_0$, i.e.:
\begin{equation}\label{pom55}
\sum_{i=t_0+1}^{m+s}{g}_i\ge \sum_{i=1}^md_i+\sum_{i=t_0+1}^s a_i,
\end{equation}
which is by the definition of ${g}_{t_0+1},\ldots,{g}_{m+s}$, equivalent to (\ref{312}).\\

Now, let $j=t_0+1,\ldots,s$. Let $x_j=\min\{i\in\{1,\ldots,h'\}| t_i=j\}$. 
 Then, by (\ref{barh2}), the condition (\ref{maj22}) becomes
$$\sum_{i=z_{x_j}+t_{x_j}+1}^{m+s} {g}_i \ge \sum_{i=z_{x_j}+1}^m d_i + \sum_{i=j+1}^s a_i,$$
which is (by the definition of ${g}_i$'s) equivalent to 
\begin{equation}\sum_{i=x_j+1}^{h'}c^i\ge \sum_{i\ge z_{x_j}+1, i\notin \Delta}d_i+\sum_{i=t_{x_j}+1}^s a_i.\label{17}\end{equation}

In order to prove (\ref{17}) we need to consider the following three possibilities:
\begin{eqnarray}
\bullet&& w_{x_j}>0, \textrm{ i.e. }  c^{x_j}>d_{z_{x_j+1}-w_{x_j}+1}>c^{x_j+1}, \textrm{ and } z_{x_j+1}-w_{x_j}+1\notin \Delta, \nonumber\\
&&\textrm{ by the part $(b)$ of the definition of the set  } \Delta\label{3s}\\
\bullet&& w_{x_j}>0, \textrm{ i.e. }  c^{x_j}>d_{z_{x_j+1}-w_{x_j}+1}>c^{x_j+1}, \textrm{ and }  z_{x_j+1}-w_{x_j}+1\notin \Delta,\nonumber\\&&\textrm{ by the part $(a)$ of the definition of the set  } \Delta,\label{2s}\\
\bullet&& w_{x_j}=0, \textrm{ i.e. there are no } i\notin \Delta, c^{x_j}>d_i>c^{x_j+1}\label{1s}
\end{eqnarray}
\vspace{0.2cm}

 \indent First consider the case (\ref{3s}). Suppose that $w_{x_j}>0$, such that   $z_{x_j+1}-w_{x_j}+1\notin \Delta,$ $ c^{x_j}>d_{z_{x_j+1}-w_{x_j}+1}>c^{x_j+1}$, satisfies the following condition (see the part $(b)$ of the definition of the set $\Delta$ and note that $q_{z_{x_j+1}-w_{x_j}+1}=t_{x_j}$):
\begin{equation}\sum_{i=x_{j}+1}^{h'}c^i\ge d_{z_{x_{j+1}}-w_{x_j}+1}+ \sum_{i> z_{x_{j+1}}-w_{x_j}+1, i\notin \Delta}d_i+\sum_{i=t_{x_{j}}+1}^s a_i.\label{19}\end{equation}
Condition (\ref{19}) is equivalent to (\ref{17}), which finishes our proof in this case.\\

 \indent Next, we consider the case (\ref{2s}). In this case we have that  $w_{x_j}>0$, and $d_{z_{x_j+1}-w_{x_j}+1}$ is among
$ \sharp\{i| a_i>c^{x_j+1}\}- s+(h'-x_j)-\sharp\{ i\notin \Delta|d_i<c^{x_j+1}\}+1$ smallest 
$e_i$'s larger than $c^{x_j+1}$ (see the part $(a)$ of the definition of the set $\Delta$), i.e. $$d_{z_{x_j+1}-w_{x_j}+1}\in\{e_{z_{x_j+1}+t_{x_j+1}}, \ldots, e_{z_{x_j+1}+m_{x_j+1}}\}.$$

Thus, in this case we have that $t_{x_j+1}\le m_{x_j+1}$. 

Let us consider the differences $m_i-t_i$ for all $i=0,\ldots,x_j+1.$ We have that $m_{x_j+1}-t_{x_j+1}\ge 0$, and $m_0-t_0=-t_0\le 0$ (because of  Lemma \ref{t0}). Thus, there exists $v:=\max\{i\in\{0,\ldots,x_{{j}}\}|m_i-t_i\le 0\}$.
Then 
$m_{v+1}-t_{v+1}\ge 0$ and $v\le x_j$,  so we have that condition $(i)$ is satisfied for $v+1$. i.e.

\begin{equation}\sum_{i=z_{v+1}+t_{v+1}}^{z_{v+1}+m_{v+1}}e_i\le\sum_{i=v+1}^{h'} c^i-\sum_{i>z_{v+1},i\notin \Delta} d_i-\sum_{i=m_{v+1}+1}^s a_i.\label{20}\end{equation}\\


Before proceeding we shall prove formulas (\ref{sept4}) and (\ref{sept5}) below:


Let $i\in\{0,\ldots,h'-1\}$.

\begin{equation}\label{sept4}\textrm{If  }\quad  m_i-t_i\le 0, \textrm{  then } c^i\ge e_{z_{i+1}+t_{i+1}}.\end{equation} 

Last is true since $z_{i+1}+t_{i+1}\ge z_i+w_i+t_i+1-w_i>z_i+m_i$.\\

On the other hand, if $m_i>t_i$, we have
$m_{i+1}-t_{i+1}+1=m_i+\sharp\{j| c^i\ge a_j> c^{i+1}\}-t_i+w_i>\sharp\{j| c^i>a_j\ge c^{i+1}\}+w_i$. Therefore $m_{i+1}\ge t_{i+1}$ and $m_{i+1}-t_{i+1}+1$ is strictly bigger than the number of $a_l$'s and $d_j$'s with $j\notin \Delta$, that are between $c^i$ and $c^{i+1}$. Therefore at least one among $e_{z_{i+1}+t_{i+1}},\ldots,e_{z_{i+1}+m_{i+1}}$ is bigger than $c^i$, i.e. $c^i<e_{z_{i+1}+t_{i+1}}$. Thus, we have 
\begin{equation}\label{sept5}\textrm{If  } \quad m_i-t_i> 0, \textrm{  then } c^i< e_{z_{i+1}+t_{i+1}}.\end{equation} \\

Now we go back to the proof of (\ref{maj22}) in the case (\ref{2s}).

First suppose that  $v=x_j$. 
Then $m_{x_j}-t_{x_j}\le 0$. This implies that $c^{x_j}\ge e_{z_{x_j+1}+t_{x_j+1}}$.
Thus, there are exactly $w_{x_j}$ of $d_i$'s among $e_{z_{x_j+1}+t_{x_j+1}}, \ldots, e_{z_{x_j+1}+m_{x_j+1}}$, and those are $d_{z_{x_j+1}-w_{x_j}+1},\ldots,d_{z_{x_j+1}}$. The remaining $m_{x_j+1}-t_{x_j+1}+1-w_{x_j}=m_{x_j+1}-t_{x_j}$ are $a_i$'s, i.e. $a_{t_{x_j}+1},\ldots,a_{m_{x_j+1}}.$ Then (\ref{20}) becomes (note that we are in the case $v=x_j$)
\begin{equation}\sum_{i=x_j+1}^{h'} c^i\ge \sum_{i>z_{x_j},i\notin \Delta} d_i+\sum_{i=t_{x_j}+1}^s a_i,\label{20a}\end{equation}
as wanted.\\






Next, suppose that   $0\le v<x_j$. 
In this case $m_i-t_i>0$, for all $i=v+1,\ldots,x_j$, and so we have that $c^i<e_{z_{i+1}+t_{i+1}}$, for all $i=v+1,\ldots,x_j$. This implies that there are no $j\in \Delta$ with $c^{v+1}>d_j>c^{x_j+1}$, and so   
$w_i=z_{i+1}-z_i$  and 
\begin{equation}\label{102}
z_{i+1}+t_{i+1}=z_i+t_{i+1}+w_i=z_i+t_i+1,\quad i=v+1,\ldots,x_j. \end{equation}
It also means that (\ref{20}) can be re-written as : 
\begin{equation}\label{101uj}
\sum_{i=z_{v+1}+t_{v+1}}^{z_{x_j}+m_{x_j}}e_i\le \sum_{i=v+1}^{h'} c^i-\sum_{i>z_{x_j}, i\notin \Delta}d_i-\sum_{i=m_{x_j}+1}^sa_i.
\end{equation}
Since $m_v-t_v\le 0$, we have 
$c^v\ge e_{z_{v+1}+t_{v+1}}$, and so
$c^v\ge e_{z_{v+1}+t_{v+1}}\ge\cdots\ge e_{z_{x_j}+m_{x_j}}>c^{x_j}$.

From the definition of $x_j$, we have $t_r<t_{x_j}=j,$ for all $r<x_j$, i.e. 
\begin{equation}\sharp\{i\notin \Delta| c^r>d_i>c^{x_j}\}< x_j-r,\quad\textrm{ for all }\quad  r<x_j.\label{s11}\end{equation}

Therefore among
$e_{z_{v+1}+t_{v+1}},\ldots, e_{z_{x_j}+m_{x_j}}$ there is at most $x_j-v-1$ $d_i$'s (note that as we have shown above, all $d_i$'s among those $e_i$'s satisfy $i\notin \Delta$). Also by (\ref{102}), $z_{x_j}+t_{x_j}=z_{v+1}+t_{v+1}+x_j-(v+1)$. Thus, among those $e_i$'s there are at  least $z_{x_j}+m_{x_j}-(z_{v+1}+t_{v+1})+1-(x_j-v-1)=z_{x_j}+m_{x_j}+1-(z_{x_j}+t_{x_j})=m_{x_j}-t_{x_j}+1$, $a_i$'s. Thus $a_{t_{x_j}}, \ldots,a_{m_{x_j}}$ surely belong to them. Since $a_{t_{x_j}}\ge a_{m_{x_j}}>c^{x_j}$ and since $e_{z_{i+1}+t_{i+1}}>c^i$, $i=v+1,\ldots,x_j$,  (\ref{s11}) and (\ref{101uj}) give

$$\sum_{i=x_j+1}^{h'} c^i\ge \sum_{i>z_{x_j},i\notin \Delta}d_i+\sum_{i=t_{x_j}+1}^sa_i,$$
i.e. we have proved (\ref{17}).\\

 \indent  So, we are left with the case (\ref{1s}), i.e. $w_{x_j}=0$, which means that  there are no  $i\notin \Delta,$ such that  $c^{x_j}>d_i>c^{x_j+1}.$

In this case, we are left with two possibilities
\begin{eqnarray}
&&t_{x_j+1}\le m_{x_j+1}\label{4s}\\
&&t_{x_j+1}> m_{x_j+1}\label{5s}
\end{eqnarray}

The case (\ref{4s}) is done exactly as in the case (\ref{2s}) when $w_{x_j}>0$ and $t_{x_j+1}\le m_{x_j+1}$.\\

So we are left with the case (\ref{5s}). The proof of this case goes by the induction on $j=t_0+1,\ldots,s$. 

Let $j=s$. Since   (\ref{novnov}) gives $c^{h'}\ge a_s$, (\ref{txs}) implies  
$t_{x}<s$ for $x<h'$. So since $t_{h'}=s$, we have $x_s=h'$. Hence (\ref{17}) becomes $0\ge 0$, which is trivially satisfied.\\

Now, fix $j\in\{t_0+1,\ldots,s-1\}$, and suppose that (\ref{17}) is satisfied for all $j+1,\ldots,s$. We shall prove that it is then also valid for $j$. 

Since $t_{x_j+1}>m_{x_j+1}$, we have  $c^{x_j+1}\ge a_{m_{x_j+1}+1}\ge a_{t_{x_j+1}}$.  Since  there are no $i\notin \Delta$ such that $c^{x_j}>d_i>c^{x_j+1}$, we have $t_{x_j+1}=t_{x_j}+1=j+1$, and so $x_{j+1}=x_j+1$. By the induction hypothesis for $j+1$, we have  
\begin{equation}\sum_{i=x_{j+1}+1}^{h'}c^i\ge \sum_{i\ge z_{x_{j+1}}+1, i\notin \Delta}d_i+\sum_{i=t_{x_{j+1}}+1}^s a_i.\label{18}\end{equation}
Since  $c^{x_j+1}\ge a_{m_{x_j+1}+1}\ge a_{t_{x_j+1}}=a_{t_{x_j}+1}$, then (\ref{18})
gives (\ref{17}).\\


%


This finishes our proof of (\ref{17}), and consequently of (\ref{maj22}).\\
\\

Finally, (\ref{maj33}) follows from (\ref{pom55}) (i.e. (\ref{312})), together with (\ref{defbg1}). Therefore we have shown that 
$$\mathbf{{g}}\prec''(\mathbf{d},\mathbf{a}).$$

\noindent Dually we obtain $$\mathbf{{g}}\prec''(\mathbf{c},\mathbf{b}).$$

This finishes the proof of Theorem \ref{glavna}.  \kraj

\section{Proof of Theorem \ref{glavna1}}
\subsection{Necessity of conditions (\ref{sumasuma}), $(i)$ and $(ii)$}\label{nesw}
Let  there exists a partition $\mathbf{g}$ such that 
\begin{equation}\mathbf{g}\prec'(\mathbf{d},\mathbf{a})\quad\textrm{ and }\mathbf{g}\prec'(\mathbf{c},\mathbf{b}).\label{est112}\end{equation}
Then (\ref{sumasuma}) follows trivially. Also, 
then $\mathbf{g}$ also satisfies 
\begin{equation}\mathbf{g}\prec''(\mathbf{d},\mathbf{a})\quad\textrm{ and }\mathbf{g}\prec''(\mathbf{c},\mathbf{b}),\label{est112uqqq}\end{equation}
and so by Theorem \ref{glavna} we obtain conditions $(i)$ and $(ii)$, as wanted.

\subsection{Sufficiency of conditions (\ref{sumasuma}), $(i)$ and $(ii)$}\label{sw}

Let  us  assume that conditions  (\ref{sumasuma}), $(i)$ and $(ii)$ are valid. By Theorem \ref{glavna} conditions $(i)$ and $(ii)$  imply the existence of a partition $\mathbf{\bar{g}}$ such that
\begin{equation}\mathbf{\bar{g}}\prec''(\mathbf{d},\mathbf{a})\quad\textrm{ and }\quad
\mathbf{\bar{g}}\prec''(\mathbf{c},\mathbf{b}).\label{est112u}\end{equation}

The rest of this section is completely analogous to  \cite{ela}. It doesn't depend on the definitions of the sets $S$ and $\Delta$, so it remains completely the same. Thus, let 
$\Omega:=  \sum_{i=1}^{m+s}\bar{g}_i -( \sum_{i=1}^s a_i + \sum_{i=1}^m d_i) \ge 0$ and let $f :=\min\{ i|\sum_{j=1}^i\bar{g}_j-i\bar{g}_i\ge \Omega\}.$ 


We are going to define $g_i$, $i = 1,\ldots,m+s,$ such that
\begin{equation}\sum_{i=1}^{m+s} g_i = \sum_{i=1}^{m}d_i + \sum_{i=1}^{s}a_i, \label{bas}\end{equation}
$$g_i=\bar{g}_i,\quad \textrm{  for all } \quad i\ge f,$$
$$\bar{g}_{f-1} \ge g_i \ge \bar{g}_f\quad\textrm{  for all }\quad i=1,\ldots,f-1,$$
 and
$$g_1 \ge g_{f-1} \ge g_1 - 1.$$
In other words, we decrease the smallest possible number of $\bar{g}_i$'s, such that the sum is correct, and such that $g_1 \ge g_2 \ge\cdots\ge g_{f-1}$ becomes the most homogeneous partition of $\bar{g}_1 + \bar{g}_2 + \cdots+ \bar{g}_f -1 - \Omega$. By Lemma \ref{2.4} such defined $g_1 \ge \cdots\ge  g_{m+s}$ satisfy 
$$\mathbf{g}\prec''(\mathbf{d},\mathbf{a})\quad\textrm{ and }\quad
\mathbf{{g}}\prec''(\mathbf{c},\mathbf{b}).$$
However, since  (\ref{sumasuma}) and (\ref{bas}) are valid, by the definition of the generalized majorization we also have
$$\mathbf{g}\prec'(\mathbf{d},\mathbf{a})\quad\textrm{ and }\quad
\mathbf{{g}}\prec'(\mathbf{c},\mathbf{b}),$$ as wanted.

\kraj

\end{document}